\documentstyle[12pt]{article}
\def\boldmat{\bf}
\def\Re{\mathop{\rm Re}\nolimits}
\def\Hom{\mathop{\rm Hom}\nolimits}
\def\End{\mathop{\rm End}\nolimits}
\def\dist{\mathop{\rm dist}\nolimits}
\def\Span{\mathop{\rm Span}\nolimits}
\def\Ker{\mathop{\rm Ker}\nolimits}
\def\Im{\mathop{\rm Im}\nolimits}
\def\tr{\mathop{\rm tr}\nolimits}
\def\ess{\mathop{\rm ess}\nolimits}
\def\Sp{\mathop{\rm Sp}\nolimits}
\def\e{{\varepsilon}}
\def\t{{\tau}}
\def\th{{\theta}}
\def\la{{\langle}}
\def\ra{{\rangle}}
\def\v{{\vert}}
\def\V{{\Vert}}
\def\p{{\perp}}
\def\arr{{\longrightarrow}}
\def\a{{\alpha}}
\def\g{{\gamma}}
\def\l{{\lambda}}
\def\i{{\infty}}
\date{22 July 1994}
\author{V.~M.~Manuilov}
\title{Diagonalization of compact operators in Hilbert
       modules over finite $W^{*}$-algebras}
\begin{document}
\maketitle
\begin{abstract}

It is known that a continuous family of compact operators can be
diagonalized pointwise. One can consider this fact as a possibility
of diagonalization of the compact operators in Hilbert modules over
a commutative $W^{*}$-algebra. The aim of the present paper is to
generalize this fact for a finite $W^{*}$-algebra $A$ not necessarily
commutative. We prove that for a compact operator $K$ acting in
the right Hilbert $A$-module $H^{*}_{A}$ dual to $H_{A}$ under slight
restrictions one can find a set of ``eigenvectors'' $x_{i}\in
 H^{*}_{A}$ and
a non-increasing sequence of ``eigenvalues'' $\l_{i}\in
 A$ such that $K\,x_{i} = x_{i}\,\l _{i}$
and the autodual Hilbert $A$-module generated by
these ``eigenvectors'' is the whole $H_{A}^{*}$. As an application we
consider the Schr\"odinger operator in magnetic field with
irrational magnetic flow as an operator acting in a Hilbert
module over the irrational rotation algebra $A_{\th}$ and discuss
the possibility of its diagonalization.
\end{abstract}

\section{Introduction}

 The classical Hilbert-Schmidt theorem states that any compact
self-adjoint operator acting in a Hilbert space can be
diagonalized. It is also known that a continuous family of compact
operators is diagonalizable. When active study of Hilbert modules
began some results were obtained concerning diagonalizability of
some operators acting in these modules. R.~V.~Kadison \cite{k1},\cite{k2}
proved that a self-adjoint operator in a free finitely generated module
over a $W^{*}$-algebra is diagonalizable. Later on some other
interesting results about diagonalization of operators appeared
\cite{gp},\cite{mur},\cite{z}. This paper is a step in the same direction
and is
concerned with diagonalization of compact operators in the Hilbert
module $H^{*}_{A}$ over a finite $W^{*}$-algebra $A$. Its main results were
announced in \cite{m2}.

\bigskip

    The present paper is organized as follows:
At section 2 we study some properties of Hilbert modules over
finite $W^{*}$-algebras related with orthogonal complementability.
The main technical result is the isomorphy of $H^{*}_{A}$ and the
orthogonal complement to $A$ in $H^{*}_{A}$.
At section 3 we recall the basic facts about the compact operators
in Hilbert modules. Here we also give an example showing that the
module $H_{A}$ is not sufficient to diagonalize compact operators,
so we must turn to its dual module $H^{*}_{A}$.
Section 4 contains the proof of the main theorem of this paper
about diagonalization of a compact operator in the module
$H^{*}_{A}$.
Here we also discuss the uniqueness condition for the
``eigenvalues'' of this operator.
Section 5 deals with quadratic forms on Hilbert modules related to
a self-adjoint operator. Properties of these forms are mostly the
same as on a Hilbert space.
At section 6 we discuss an example which motivated the present
paper. We consider the perturbated Schr\"odinger operator with
irrational magnetic flow as an operator acting in a Hilbert module
over the irrational rotation algebra $A_{\th}$ and we show that this
operator is diagonalizable.

\bigskip

{\bf Acknoledgement.\/}
This work was partially
supported by the Russian Foundation for Fundamental Research (grant
\mbox{N 94-01-00108-a)} and the International Science Foundation
(grant N MGM000).
I am indebted to M.~Frank, A.~A.~Irmatov, A.~S.~Mishchenko
and E.~V.~Troitsky for helpful discussions.

\section{Orthogonal complements in Hilbert modules over
                    finite $W^{*}$-algebras}

 Throughout this paper $A$ is a finite $W^{*}$-algebra admitting the
central decomposition into a direct integral over a compact Borel
space. By $\t$ we denote a normal faithful finite trace
on $A$ with $\t (1) = 1$. Recall some facts about Hilbert modules.
Standart references on them are \cite{ka},\cite{mi},\cite{p1}. If $B$
is a $C^{*}$-algebra
we denote by $H_{B}$ (another usual denotation is $l_{2} (B)$) the
right Hilbert $B$-module consisting of the sequences $(x_{i})$,
$\ i \in {\boldmat N}$
for which the series $\sum_{i} x_{i}^{*} x_{i}$
converges in the norm topology in $B$
with the inner product $\la x,y \ra = \sum_{i} x_{i}^{*}
 y_{i}$ and the norm $\V x \V = \V \la x,x \ra
\V ^{1/2}$. Let $H^{*}_{B}$ be its dual module,
$H^{*}_{B} = \Hom _{B}(H_{B};B)$.
It is shown in \cite{p2} that in the case of $W^{*}$-algebras the inner
product on the module $H_{B}$ can be prolonged to the inner product
on the module $H^{*}_{B}$ and this module is autodual, i.e.
$(H^{*}_{B})^{*} = H^{*}_{B}$.

    Let $M\subset H^{*}_{B}$ be a Hilbert $B$-submodule. By
$M^{\p}$ we denote
its orthogonal complement in $H^{*}_{B}$. It is well-known \cite{df} that if
$M$ is a finitely generated projective Hilbert $B$-submodule in
$H^{*}_{B}$
then it is orthogonally complementary: $H^{*}_{B} = M\oplus
 M^{\p}$. If we
change $H^{*}_{B}$ by $H_{B}$ then the orthogonal complement to $M$ in
$H_{B}$ is isomorphic to $H_{B}$, but nothing is known in general about
isomorphy between $M^{\p}$ and $H^{*}_{B}$. The following theorem solves
this problem in the case of modules over a $W^{*}$-algebra decomposable
into a direct integral of finite factors and having a faithful
finite trace.

\bigskip

{\bf Theorem 2.1.\/} {\it \/ If $M$ is a finitely generated projective
$A$-submodule in $H^{*}_{A}$ then $M^{\p}$ is isomorphic to
$H^{*}_{A}$.}

\bigskip

{\bf Proof.\/} The idea of the following proof is contained in \cite{df}. Let
$g_{1},\dots , g_{n}$ be generators of the module $M$. Without loss of
generality we can assume that the operators $\la g_{i}, g_{i}
\ra \in A$ are projections, $\la g_{i}, g_{i}\ra =
 p_{i}$.
 Let $\{ e_{m}\}$ be the standart basis of
the module $H_{A} \subset H_{A}^{*}$. Fix $\e < 0$ and define elements
$e'_{m} \in M^{\p}$ by the equality
\[
e'_{m} = e_{m} - \sum_{i=1}^{n} g_{i} \la g_{i}, e_{m}\ra.
\]
Then we have
\[
\la e'_{m}, e'_{m}\ra = 1 - \sum_{i=1}^{n} \la g_{i},
e_{m}\ra^{*} \la g_{i}, e_{m}\ra.
\]
It follows from the equality
\[
  \t\Bigl(\la g_{i},g_{i}\ra\Bigr) =
  \t\biggl(\sum_{m} \la g_{i}, e_{m}\ra^{*} \la g_{i},
  e_{m}\ra\biggr)
\]
that the series $\sum_{m}\t\Bigl(\la g_{i}, e_{m}\ra^{*}
\la g_{i}, e_{m}\ra\Bigr)$ converges and there
exists such number $m_{0}$ that for any $m > m_{0}$ the inequalities
\[
  \t\Bigl(\la g_{i}, e_{m}\ra^{*}\la g_{i}, e_{m}
  \ra\Bigr) < \frac{\e}{2n};
\]
\[
  \t\Bigl(\la e'_{m},  e'_{m}\ra\Bigr) > 1  -  \frac{
  \e}{2}
\]
hold.

\bigskip

{\bf Lemma 2.2.\/}{\it \/ If $x\in H^{*}_{A},\ \V x\V = 1$ and
$\t\left(\la x, x\ra\right) > 1 - \frac{\e}{2}$ then
there exists a projection $p\in A$ with $\t (p) > 1 - \e$ such that
$p\la x, x\ra p$ is an invertible operator in the $W^{*}$-algebra
$pAp$.}

\bigskip

{\bf Proof.\/} Let $dP(\l )$ denote the projection-valued measure for the
operator $a = \la x, x\ra\in A$; $a = \int^{1}_{0} \l
dP(\l )$. Put
\[
  f(\l ) = \left\lbrace
                \begin{array}{cc}
                0,& \l\leq\l_{0}, \\
                1,& \l > \l_{0},
                \end{array}
                \right.
\]
where $\l_{0}\in [0;1]$. Then $f(a) = p$ is a projection. Denote
$d\t\left(P(\l )\right)$ by $d\mu (\l )$. It is a usual measure on
$[0;1]$ and by \cite{mn}
\[
  \t (a) = \int^{1}_{0} \l\, d\mu (\l ).
\]
We have
\begin{equation} \label{1}
1 - \frac{\e}{2} < \t (a) = \int^{1}_{0} \l\, d\mu
(\l )
 \leq \l_{0}\,\mu\left([0;\l_{0})\right) + \mu\left(
[\l_{0};1]\right).
\end{equation}
Since $P(1) = 1$ we have
\begin{equation} \label{2}
\mu\left([0;\l_{0})\right) + \mu\left([\l_{0};1]\right) = 1
\end{equation}
{}From (\ref{1}) and (\ref{2}) we obtain the inequality
\[
\mu\left([\l_{0};1]\right) > 1 - \frac{\e}{2\left( 1 -
\l_{0}\right)}.
\]
Choosing an appropriate number $\l_{0}\neq 0$ we obtain
$\mu\left([\l_{0};1]\right) > 1 - \e$.
{}From the definition of the function $f(\l )$ we have
\begin{eqnarray*}
\t (p) & = & \t\left(f(\l )\right) \quad = \quad \int^{1}_{0}
f(\l )\,d\mu (\l ) \\
& = & \int^{1}_{\l_{0}} d\mu (\l ) \quad = \quad
\mu\left([\l_{0};1]\right)\quad > \quad 1 - \e.
\end{eqnarray*}
Consider now the operator $pap\in pAp$. The equality
\[
pap = \int^{1}_{0} \l\, f(\l )\, dP(\l )
\]
follows from the spectral theorem, therefore the spectrum of the
operator $pap$ as an element of the $W^{*}$-algebra $pAp$ lies in
$[\l_{0};1]$, hence is separated from zero and this operator
is invertible in $pAp$.$\quad\bullet$ \\
Let
\[
\Bigl(p\la e'_{m}, e'_{m}\ra p\Bigr)^{-1/2} = pbp = b\in pAp
\]
and $e''_{m} = e'_{m}\cdot b$. Then
\[
\la e''_{m}, e''_{m}\ra = pbp\la e'_{m}, e'_{m}\ra pbp
= p.
\]
Now take an element $y\in M^{\p}$, $y\neq 0$, $\V y\V\leq 1$.
For every $\e > 0$
beginning from a certain number $m$ we have
\begin{eqnarray*}
\t\Bigl(\la e''_{m},y\ra^{*}\la e''_{m},y\ra\Bigr)
& = & \t\left(\la e_{m},y\ra^{*} b^{2}\la e_{m},y\ra
\right)\\
& \leq & \V b\V^{2}\t\Bigl(\la e_{m},y\ra^{*}
\la e_{m},y\ra\Bigr) \quad < \quad \frac{\e^{2}}{2}
\end{eqnarray*}
because of convergence of the series $\sum_{m}\t\Bigl(\la
e_{m},y\ra^{*}\la e_{m},y\ra\Bigr)$.

Denote the operator $\la e''_{m},y\ra^{*}
\la e''_{m},y\ra\in A$ by $c$, then $\t (c) < \frac{\e^{2}}
{2}$;$\ \V c\V\leq 1$
 and $c\geq 0$ (i.e. $c$ is a positive operator). If $dQ(\l )$
is its projection-valued measure and if we denote the measure
$d\t (Q(\l ))$ by $d\nu (\l )$  then we have
$\int_{0}^{1} \l\, d\nu (\l ) < \frac{\e^{2}}{2}$.
If $\l_{1}\in [0;1]$ then
\[
  \l_{1}\cdot\int_{\l_{1}}^{1} d\nu (\l ) \leq
\int^{1}_{0} \l\, d\nu (\l ) = \t (c) < \frac{\e^{2}}{2},
\]
hence $\int^{1}_{\l_{1}} d\nu (\l ) < \frac{\e^{3}}{2}$.
 Taking $\l_{1} = \frac{\e^{2}}{2}$ we obtain
$\nu ([\l_{1};1]) < \e $.
 Put
\[
g(\l ) = \left\lbrace
              \begin{array}{cc}
              1, & \l < \l_{1}, \\
              0, & \l \geq \l_{1}.
              \end{array}
              \right.
\]
Then $q = g(c)$ is a projection with $\t (q) = \nu ([0;\l_{1}))
> 1 - \e $
and
\begin{equation} \label{3}
\V qcq\V \leq \l_{1} = \frac{\e^{2}}{2}.
\end{equation}
By $p\lor q$ (resp. $p\land q$) we denote the least upper (resp.
greatest lower) bound for projections $p$ and $q$. Put $p' = p\land q$.
As by \cite{t}
\[
\t (p) + \t (q) = \t (p\lor q) + \t (p\land q),
\]
so we have
\[
\t (p') = \t (p) + \t (q) - \t (p\lor q) >
(1 - \e ) + (1 - \e ) - 1 = 1 - 2\e
\]
because $\t (p\lor q) \leq \t (1) = 1$. The inequality
\begin{equation} \label{4}
\V p'cp'\V\leq\frac{\e^{2}}{2}
\end{equation}
follows from (\ref{3}). Put now $e'''_{m} = e''_{m}\cdot p'$. Then
$\la e'''_{m},e'''_{m}\ra = p'$.
Put further $y' = y + \e e'''_{m}\in M^{\p}$. We can decompose
$y'$into
two orthogonal summands: $y' = u + v$ , where $u=y-e'''_{m}
\la e'''_{m},y\ra$,$v=e'''_{m}\Bigl(\la e'''_{m},y\ra +
\e\cdot 1\Bigr)$; $u,v\in M^{\p}$. Then
\[
\la y',y'\ra = \la u,u\ra + \Bigl(\la e'''_{m},y\ra
+ \e p'\Bigr)^{*}\Bigl(\la e'''_{m},y\ra + \e p'\Bigr)
\]
and
\begin{eqnarray} \label{5}
p'\la y',y'\ra p' & = & p'\la u,u\ra p' +
\Bigl(p'\la e'''_{m},y\ra p' + \e p'\Bigr)^{*}\Bigl(p'\la e'''_{m},
y\ra p' + \e p'\Bigr) \nonumber \\
& = & p'\la u,u\ra p' + \Bigl(\la e'''_{m},y\ra p' +
\e p'\Bigr)^{*}\Bigl(\la e'''_{m},y\ra p' + \e p'\Bigr).
\end{eqnarray}
Since
\begin{eqnarray*}
\Bigl(\la e'''_{m},y\ra p'\Bigr)^{*} \la e'''_{m},y\ra p' & = &
p'\la e''_{m},y \ra^{*}p'\la e''_{m},y\ra p' \\
& \leq & p'\la e''_{m},y\ra^{*}\la e''_{m},y\ra p'
\quad = \quad p'cp'
\end{eqnarray*}
it follows from (\ref{4}) that $\V \la e''_{m},y\ra p'\V
\leq \frac{\e}{\sqrt{2}} < \e$. Therefore
the operator $\la e''_{m},y\ra p' + \e p'$ is invertible
in the $W^{*}$-algebra $p'Ap'$.
The invertibility of $p'\la y',y'\ra p'$ follows now from (\ref{5}).
Consider the trace norm on $M^{\p}$ (and on $H^{*}_{A}$) defined by
\[
\V x\V_{\t} = \t \Bigl(\la x,x\ra\Bigr)^{1/2}.
\]
The inequality
\[
\t \Bigl(\la y' -y,y'-y\ra \Bigr) = \t \Bigl(\e^{2}\la
e'''_{m},e'''_{m}\ra \Bigr) = \t (\e^{2} p') < \e^{2}
\]
gives us the estimate $\V y'-y\V_{\t} < \e$.
So we have proved that the
elements of $M^{\p}$     for which there exists a projection $p'$
 with $\t (p') > \frac{1}{2}$
such that $p'\la x,x\ra p'$ is invertible in $p'Ap'$ are dense in
$M^{\p}$ in the trace norm.

\bigskip

{\bf Corollary 2.3.\/}$\,${\it There exists some $x\in M^{\p}$ such that
$\V x\V_{\t} > \frac{1}{2}$ and $\V x\V\leq 1$.}

\bigskip

Let now $\lbrace y_{n}\rbrace$ be a sequence containing every $e_{m}$
infinitely
many times. Put $y=y_{1}-\sum_{k=1}^{n} g_{k}\la g_{k},y_{1}
\ra$. Then for $\e_{1}=1$ there
exists some $y'\in M^{\p}$  with $\V y'\V\leq 1$ such that
\begin{equation} \label{6}
\V y-y'\V_{\t} < \e_{1}
\end{equation}
and a projection $p_{1}$ with $\t (p_{1})>\frac{1}{2}$ such that
$p_{1}\la y',y'\ra p_{1}$ is
invertible in $p_{1}Ap_{1}$. Then putting $h_{1}=y'b'$ where
$b'=\Bigl(p_{1}\la y',y'\ra p_{1}\Bigr)^{-1/2}\in p_{1}Ap_{1}$
                           we obtain from (\ref{6}) the inequality
\[
\dist_{\t }(y,B_{1}(h_{1}A))\leq\dist_{\t }(y,h_{1}(b')^{-1})=
\dist_{\t }(y,y')<\e_{1},
\]
where by $B_{1}$ we denote the unit ball of a Hilbert module in the
initial norm. Therefore
\[
\dist_{\t }(y_{1},B_{1}(\Span_{A}(M,h_{1})))<\e_{1}.
\]
Then taking $\e_{2}=\frac{1}{2}$ we can find an element
$h_{2}\in (\Span_{A}(M,h_{1}))^{\p}$
such that $\la h_{2},h_{2}\ra=p_{2}$ is a projection with
$\t (p_{2})>\frac{1}{2}$ and
\[
\dist_{\t }(y_{2},B_{1}(\Span_{A}(M,h_{1},h_{2})))<\e_{2}.
\]
Continuing this process and taking $\e_{k}=\frac{1}{k}$
we obtain a set of
mutually orthogonal elements $h_{i}\in M^{\p}$ with
$\la h_{i},h_{i}\ra=p_{i}$
being a projection and $\t (p_{i})>\frac{1}{2}$ such that
\begin{equation} \label{7}
\dist_{\t }(y_{k},B_{1}(\Span_{A}(M,h_{1},h_{2},\ldots ,h_{k})))
<\e_{k}
\end{equation}
These $h_{i}$ generate an $A$-module $N\subset M^{\p}$
and from (\ref{7}) we have
\[
\dist_{\t }(y_{k},B_{1}(M\oplus N))<\frac{1}{k},
\]
hence the trace norm closure of $B_{1}(M\oplus N)$ contains the unit
ball of the whole $H^{*}_{A}$ and the trace norm closure of $B_{1}(N)$
contains $B_{1}(M^{\p})$.

    The constructed above basis $\lbrace h_{i}\rbrace$ of $N$ is inconvenient
because the inner squares of $h_{i}$ are not unities. So we have to
alter it. By $T$ we denote the standart center-valued trace on $A$.

\bigskip

{\bf Lemma 2.4.\/}{\it \/ For any number $C$ there exists some number
$n$ such that $T\left(\sum^{n}_{i=1}p_{i}\right)\geq C$.}

\bigskip

{\bf Proof.\/} Suppose that there exists a normal state $f$ on the center
$Z$ of $A$ such that for some $C\quad (f\circ T)\left(\sum_{i=1}^{\i}
p_{i}\right)<C$.
 Then there exists
a central projection $z\in Z$ such that
\begin{equation} \label{8}
T\biggl(\sum_{i=1}^{\i}p_{i}z\biggr)<C.
\end{equation}
Consider the $W^{*}$-algebra $zAz$ . Multiplication by $z$ turns any
Hilbert module over $A$ into a Hilbert module over $zAz$ and
preserves orthogonality of submodules. So we have $Nz\subset M^{\p}z$
and $B_{1}(Nz)$ is dense in $B_{1}(M^{\p}z)$ in the trace norm
$\V\cdot\V_{\t_{z}}$
defined by the faithful trace $\t_{z}$ on $zAz$ induced by $\t$ . The
inequality (\ref{8}) means that for any $\e >0$ changing $z$ by a lesser
central projection if nessessary we can find such number $k$ that
the inequality $T\left(\sum_{i>k}p_{i}z\right)<\e$
holds. Decompose the module $Nz$ :
$Nz=L_{k}\oplus R_{k}$ where $L_{k}$ is the $zAz$-module generated by
$h_{1}z,\ldots ,h_{k}z$ and $R_{k}$ is the orthogonal complement to $L_{k}$
in
$Nz$. As $B_{1}(Nz)$ is dense in $B_{1}(M^{\p}z)$ so $B_{1}(R_{k})$ must be
dense in $B_{1}((Mz\oplus L_{k})^{\p})$ in the trace norm $\V\cdot
\V_{\t_{z}}$. Let $x=\sum_{i>k}h_{i}x_{i}\in B_{1}(R_{k})$.
 Estimate its trace norm:
\begin{eqnarray*}
\V x\V^{2}_{\t_{z}} & = & \t_{z}\biggl(\sum_{i>k}x^{*}_{i}
\la h_{i}z,h_{i}z\ra x_{i}\biggr)\quad = \quad
\t_{z}\biggl(\sum_{i>k}x^{*}_{i}p_{i}zx_{i}\biggr) \\
& = & \!\!\t_{z}\biggl(\sum_{i>k}p_{i}zx_{i}x^{*}_{i}\biggr)\quad\leq\quad
\t_{z}\biggl(\sum_{i>k}p_{i}z\cdot\V x_{i}\V^{2}\biggr) \\
& \leq & \!\!\t_{z}\biggl(\sum_{i>k}p_{i}z\cdot\V x\V^{2}\biggr)
\quad \leq\quad \t_{z}\biggl(\sum_{i>k}p_{i}z\biggr).
\end{eqnarray*}
As we have $T\left(\sum_{i>k}p_{i}z\right)<\e$ so
$\t_{z}\left(\sum_{i>k}p_{i}z\right)<\e$ and so we
obtain $\V x\V^{2}_{\t_{z}}<\e$ for all $x\in B_{1}(R_{k})$.
But as $B_{1}(R_{k})$
is dense in $B_{1}((Mz\oplus L_{k})^{\p})$ so for all $y\in
B_{1}((Mz\oplus L_{k})^{\p})$
we have $\V y\V^{2}_{\t_{z}}<\e$. On the other
hand if we apply the
corollary 2.3 to the module $(Mz\oplus L_{k})^{\p}$ instead of
$M^{\p}$ we can
find in $B_{1}((Mz\oplus L_{k})^{\p})$ an element $y$ with
$\V y\V_{\t_{z}}>\frac{1}{2}$. The
obtained contradiction finishes the proof.$\quad \bullet$

    Choose now a projection $q$ in $A$ with the properties:
\begin{equation} \label{9}
T(q)=\min (T(p_{1}+p_{2});1) - T(p_{1})
\end{equation}
and $q\p p_{1}$. It follows from (\ref{9}) that $T(q)\leq T(p_{2})$, therefore
there exists another projection $q'$ equivalent to $q$ such that $q'\leq
p_{2}$.
Equivalence of $q$ and $q'$ involves existance of a unitary $u\in A$
such that $qu=uq'$. Put $r=h_{2}q'u^{*}\in N$. Then $r$ is
orthogonal to $h_{1}$ and
\[
\la r,r\ra =uq'\la h_{2},h_{2}\ra q'u^{*}=uq'p_{2}q'u^{*}
=uq'u^{*}=q.
\]
Put further $H_{1}^{(1)}=h_{1} + r$. Then
\[
\la h_{1}^{(1)},h_{1}^{(1)}\ra =\la h_{1},h_{1}\ra +
\la r,r\ra =p_{1} +q.
\]
Notice that $T(\la h_{1}^{(1)},h_{1}^{(1)}\ra ) =
\min (T(p_{1} +p_{2});1)$. Taking into
consideration the next element $h_{3}$ we can obtain $h_{1}^{(2)}$ such that
$\la h_{1}^{(2)},h_{1}^{(2)}\ra $ is a projection and
$T(\la h_{1}^{(2)},h_{1}^{(2)}\ra ) = \min (T(p_{1} +p_{2}
+p_{3});1)$.
Repeating this procedure and increasing the value of
$T(\la h_{1}^{(n)},h_{1}^{(n)}\ra ) $
we can construct by the lemma 2.4 an element $h_{1}^{\i}$ such that
$\la h_{1}^{\i},h_{1}^{\i}\ra =1$.
The orthogonal complement to $h_{1}^{\i}$ in $N$ is generated by
elements $h_{i}q_{i},\ i>1$ where $q_{i}$ are some projections. Applying
the construction described above to these generators we can
construct by induction a set of elements $h_{i}^{\i}$ with
$\la h_{i}^{\i},h_{i}^{\i}\ra =1$
which generates the module $N$. Hence $\{ h_{i}^{\i}\}$ is a basis in
$N$ and $N$ is isomorphic to $H_{A}$.

    Finally we must prove that $N^{*}=M^{\p}$. As $N\subset M^{\p}$
is closed
in the usual norm, so for any $f\in (M^{\p})^{*}$ its restriction
$f\v_{N}$ belongs to $N^{*}$. Notice that the module $M^{\p}$
is autodual, $(M^{\p})^{*} = M^{\p}$
because of autoduality of $H^{*}_{A}$ and $M$. Suppose that $f\v_{N}=0$.
Since $N$ is dense in $M^{\p}$ in the trace norm, we have $f=0$ on
$M^{\p}$
because of continuity of the map $f:M^{\p}\arr A$ in this norm due
to the inequality
\begin{equation} \label{10}
\t \Bigl((f(y))^{*}f(y)\Bigr)\leq \V f\V^{2}\cdot\t \Bigl(\la y,y\ra \Bigr)
\end{equation}
where $y\in M^{\p}$. So monomorphity of the map
$M^{\p}\arr N^{*}$ is proved.
Let now $\phi\in N^{*}$. This functional can be prolonged to a map from
$M^{\p}$ to $A$. If $\{y_{n}\}\subset N$ is a sequence converging to
$y\in M^{\p}$
in the trace norm then put $\phi (y)=\lim\phi (y_{n})$. Correctness of
this definition follows from (\ref{10}) with $\phi$ instead of $f$. So the
$A$-modules $M^{\p}$ and $N^{*}$ coinside and the theorem is proved because
the module $N^{*}$ is isomorphic to $H^{*}_{A}$.$\quad\bullet$

\bigskip

{\bf Proposition 2.5.\/}{\it \/ Let $N\subset H^{*}_{B}$ be a Hilbert
submodule over a
$W^{*}$-algebra $B$ and let $N^{\p}=0$. Then its dual module $N^{*}$
coinsides with $H^{*}_{B}$.}

\bigskip

{\bf Proof.\/} According to supposition for any $z\in H^{*}_{B}$ there exists
some
$x\in N$ such that $\la z,x\ra\neq 0$. Therefore the map
$z\longmapsto \la z,\cdot\ra$
defines the monomorphism $j^{*}:H^{*}_{B}\arr N^{*}$ which is dual to the
inclusion $j:N\hookrightarrow H^{*}_{B}$ after identification of $H^{*}_{B}$
and its
dual $(H^{*})^{*}$. Their composition
\[
j^{*}\circ j: N\arr H^{*}_{B}\arr N^{*}
\]
coinsides with the natural inclusion $N\hookrightarrow N^{*}$. Its dual map
\[
(j^{*}\circ j)^{*}=j^{*}\circ j^{**}:N^{*}=N^{**}\arr
H^{*}_{B}\arr N^{*}
\]
must be an isomorphism, therefore $j^{*}$ must be epimorphic.$\quad\bullet$

\bigskip

{\bf Proposition 2.6.\/}{\it \/ Let $B$ be a $W^{*}$-algebra and let
$R\subset H_{B}$ be a
$B$-submodule without  orthogonal  complement,  i.e.  $R^{\p}=0$ in
$H_{B}$.
Then $R^{*}=H^{*}_{B}$. }

\bigskip

{\bf Proof.\/} It is easy to verify that if $R\subset H_{B}$ then
$R^{*}\subset H^{*}_{B}$. As
the module $R^{*}$ is autodual, so by \cite{f1} $R^{*}$ is orthogonally
complementary, therefore $H_{B}=R^{*}\oplus S$ with some $B$-module $S$.
Notice that the map $H_{B}\arr R^{*}$;$\ x\longmapsto\la
x,\cdot\ra$ is monomorphic by
supposition. So we have $S\p H_{B}$. But as it is known that
$H^{\p}_{B}=0$
in $H^{*}_{B}$, so $S=0$.$\quad\bullet$

\section{Compact self-adjoint operators in Hilbert \mbox{ $A$-modules}}

\setcounter{equation}{0}

 By $\End^{*}_{B}(M)$ we denote the set of all bounded $B$-linear
operators acting on a Hilbert $B$-module $M$ over a $C^{*}$-algebra $B$
and possessing a bounded adjoint operator.

\bigskip

{\bf Proposition 3.1.\/}{\it \/ If $B$ is a $W^{*}$-algebra then
$\End^{*}_{B}(H^{*}_{B})$ is a $W^{*}$-algebra.}

\bigskip

{\bf Proof\/} is reduced to verification of the isomorphy between
$\End^{*}_{B}(H^{*}_{B})$
and the $W^{*}$-algebraic tensor product of $B$ by the algebra of
bounded operators on the separable Hilbert space.

    Recall the definition of the compact operators in a Hilbert
$B$-module $M$. Put $\th_{x,y}(z)=x\la y,z\ra$ for $x,y,z\in M$.
Then $\th_{x,y}\in\End^{*}_{B}(M)$. The set ${\boldmat K}(M)$
of compact operators is
the norm-closed linear hull of the set of all operators of the form
$\th_{x,y}$. Denote by $L_{n}(B)$ the Hilbert $B$-submodule of the
modules $H_{B}$ or $H^{*}_{B}$ generated by the first $n$ elements of the
standart basis $e_{1},\ldots ,e_{n}$.

\bigskip

{\bf Proposition 3.2.\/}{\it \/ Let $C^{*}$-algebra $B$ be unital.
Then an operator
$K\in\End^{*}_{B}(H_{B})$ is compact if and only if the norm of the
restriction of $K$ to the orthogonal complement to $L_{n}$ tends
to zero.}

\bigskip

{\bf Proof.\/} Denote by $P_{n}$ the projection $H_{B}\arr
L_{n}(B)^{\p}$. Then for any $z\p L_{n}(B)$
               we have
\begin{eqnarray*}
\V\th_{x,y}(z)\V^{2} & = & \V\la \th_{x,y}(z),
\th_{x,y}(z)\ra\V\ =\ \V\la y,z\ra^{*}
\la x,x\ra\la y,z\ra\V \\
& \leq & \V x\V^{2}\V\la y,z\ra\V^{2}\,=\,
\V x\V^{2}\V\la P_{n}y,z\ra\V^{2} \\
& \leq &
\V x\V^{2}\cdot\V P_{n}y\V^{2}\cdot\V z\V^{2}.
\end{eqnarray*}
As $\V P_{n}y\V$ tends to zero, so does the norm of $\th_{x,y}$ restricted
to $L_{n}(B)^{\p}$. The same is true for linear combinations of such
operators and for their norm closure. Suppose now that for some
operator $K$ we have $\V K\v_{L_{n}(B)^{\p}}\V\rightarrow 0$.
Then as $\sum_{m=1}^{n}Ke_{m}\la e_{m},z\ra =0$
    for any $z\p L_{n}(B)$, so if $\V z\V\leq 1$ and $z\p L_{n}(B)$ then
we have
\begin{equation} \label{11}
\sup\limits_{z}\V Kz-\sum_{m=1}^{n}Ke_{m}\la e_{m},z\ra\V =
\sup\limits_{z}\V Kz\V\arr 0
\end{equation}
when $n\arr\i$ . If $z\in L_(B)$ then $Kz=
\sum_{m=1}^{n}Ke_{m}\la e_{m},z\ra$. It
means that (\ref{11}) holds also if the supremum is taken in the unit
ball of the whole $H_{B}$, therefore the operator $K$ is the norm
topology limit of the operators $K_{n}=\sum_{m=1}^{n}\th_{Ke_{m},e_{m}}.
\nolinebreak\cite{df}\quad\bullet$

\bigskip

{\bf Remark 3.3.\/} This property of the compact operators was taken as
their definition in \cite{mf}. Without the supposition that $B$ is
unital these two definitions fail to be equivalent. As it was shown
in \cite{f2} the property of an operator to be compact strongly depends
on the choice of a Hilbert structure. Throughout this paper we
consider only the standart Hilbert structure on $H_{B}$.

\bigskip

Let $K$ be a self-adjoint compact operator acting in $H_{A}$. Due
to its self-adjointness this operator can be prolonged to an
operator $K^{*}$ in $H^{*}_{A}$.

\bigskip

{\bf Lemma 3.4.\/}{\it \/ If $\,\Ker K=0$ in $H_{A}$ then $\Ker K^{*}=0$ in
$H^{*}_{A}$.}

\bigskip

{\bf Proof\/} obviously follows from the proposition 2.6. One must take the
norm closure of $\Im K$ in $H_{A}$ as $R$. Then $\Ker K=R^{\p}=0$,
hence $R^{*}=H^{*}_{A}$ and $\Ker K^{*}=(R^{*})^{\p}=0.\quad\bullet$

\bigskip

For now on we shall not distinguish the operator $K$ and its
prolongation $K^{*}$ and denote both of them by $K$.

\bigskip

Now we shall produce an example which shows the necessity of
consideration of the dual Hilbert modules if we want to diagonalize
compact operators.

\bigskip

{\bf Example 3.5.\/} Let $A=L^{\i}([0;1])$ and let $b_{k}$ be a monotonous
sequence of positive numbers converging to zero. Put
\[
a_{k}=\left\{
             \begin{array}{cc}
             1, & t\in (\frac{1}{2^{k}};\frac{1}{2^{k-1}}], \\
             0, & \mbox{for other} \quad t,
             \end{array}
      \right.
\]
and put $f_{k}(t)=b_{k}\cdot a_{k}(t)$. Let $K$ be a compact operator
which can be written in the form
\[
K=\left(
        \begin{array}{ccccc}
        f_{1} & f_{2} & \cdots & f_{n} & \cdots \\
        f_{2} & 0 & \cdots & 0 & \cdots \\
        \vdots & \vdots &  & \vdots & \\
        f_{n} & 0 & \cdots & 0 & \cdots \\
        \vdots & \vdots &  & \vdots &
        \end{array}
  \right)
\]
in the standart basis of $H_{A}$. One can easily diagonalize
pointwise this operator. Then the eigenvector corresponding to the
maximal eigenvalue can be written as $x=(x_{n}(t))$ with
$x_{1}(t)=a_{1}(t) + \frac{\sqrt{2}}{2}\sum_{k>1}a_{k}(t)$,
 and $x_{n}=\frac{\sqrt{2}}{2} a_{n}(t)$ when $n>1$. Then
$\la x,x\ra =\sum_{k} a_{k}(t) =1$.
This series converges but not in the norm topology of
$A$, so we have $x\in H^{*}_{A} \setminus H_{A}$.

\section{Diagonalization of compact operators in $H^{*}_{A}$}

\setcounter{equation}{0}

We say that a compact operator $K$ in a Hilbert module $M$ is positive if for
any $x\in M$ the operator $\la Kx,x\ra\in A$ lies in the positive cone
of $A$. In Hilbert modules as well as in Hilbert spaces positive
operators are self-adjoint. A set of elements $\{x_{i}\} \in H^{*}_{A}$ we
call a {\it ``basis''}\/ if $\la x_{i},x_{i}\ra =\delta_{ij}$ and if the dual
$A$-module for the
module generated by this set coinsides with $H^{*}_{A}$, i.e.
$(\Span_{A}\{x_{i}\})^{*}=H^{*}_{A}$. Notice that a ``basis'' is neither
algebraic
nor topological basis. An element $x\in H^{*}_{A}$ we call an
{\it ``eigenvector''}\/ and an operator $\l \in A$ we call an
{\it ``eigenvalue''}\/ for
$K$ if $x$ generates a projective $A$-module and $Kx=x\l$.

\bigskip

{\bf Theorem 4.1.\/}{\it \/ Let $K$ be a compact positive operator in
$H^{*}_{A}$ with $\Ker K=0$.
Then there exists a ``basis'' $\{x_{i}\}$ in $H^{*}_{A}$ consisting
of ``eigenvectors'', i.e. $Kx_{i}=x_{i}\l_{i}$ for some ``eigenvalues''
$\l_{i}\in A$.}

\bigskip

{\bf Proof.\/} The $W^{*}$-algebra $\End^{*}_{A}(H^{*}_{A})$ is semifinite
and its center
is the same as the center $Z$ of $A$, so this algebra as well as $A$
can be decomposed into a direct integral of factors over the
compact Borel space $\Gamma$ with the finite measure $d\g$ such that
$L^{\i}(\Gamma)=Z$.
        The operator $K$  then also can be decomposed,
\[
K=\int_{\Gamma}^{\oplus}K(\g)\,d\g .
\]
If we put $\bar T =T\otimes \tr$ where $T$ is the standart $Z$-valued finite
trace on $A$ and $\tr$ is the standart trace in the Hilbert space we
obtain a semifinite center-valued trace on the $W^{*}$-algebra
$\End^{*}_{A}(H^{*}_{A})$.
At first we show that if we separate the spectrum of $K$ from zero
then we find ourselves in the finite trace ideal of
$\End^{*}_{A}(H^{*}_{A})$.
Let $\chi_{E}$ denote as usual the characteristic function of a set
$E\subset {\boldmat R}$.

\bigskip

{\bf Lemma 4.2.\/}{\it \/ For every $\e >0$ almost everywhere on
$\Gamma$ we have $\bar T (\chi_{(\e ;+\i)}(K))<\i$.}

\bigskip

{\bf Proof.\/} Denote the spectral projection $\chi_{(\e ;+\i)}(K)$ by
$P$. Then the operator inequality
\begin{equation} \label{12}
K\v_{\Im P}\geq\e
\end{equation}
is satisfied on the $A$-submodule $\Im P\subset H^{*}_{A}$ by the spectral
theorem. Due to compactness of $K$ we can decompose $H_{A}$ into a
direct sum: $H_{A}=L_{n}(A)\oplus R$ with such number $n$ that
$\V K\v_{R}\V <\e$.
If we pass on to the dual modules then we obtain the estimate
\begin{equation} \label{13}
\V K\v_{R^{*}}\V <\e
\end{equation}
where $H^{*}_{A}=L_{n}(A)\oplus R^{*}$. Denote by $Q$ the projection in
$H^{*}_{A}$ onto $R^{*}$. Then the projection onto $\Im P\cap R^{*}$
will be $P\land Q$ and
the projection onto $(\Ker P\cap L_{n}(A))^{\p}$ will be $P\lor Q$. By the
results of \cite{t} we have
\begin{equation} \label{14}
\bar T (P\lor Q)=\bar T (P-P\land Q).
\end{equation}
As $P\lor Q\leq 1$  where $1$ stands for the unity operator in
$\End^{*}_{A}(H^{*}_{A})$,
so we obtain the inequality $\bar T (P\lor Q-Q)\leq\bar T (1-Q)$. But $1-Q$
is the projection onto $L_{n}(A)$ and its trace is equal to $n$, so from
(\ref{14}) we have $\bar T (P-P\land Q)\leq n$. Comparing (\ref{12}) with
(\ref{13}) we
conclude that $\Im P\cap R^{*}=0$, so $P\land Q=0$ and finally we have
$\bar T (P)\leq n.\quad\bullet$

\bigskip

We shall need subsequently one simple fact concerning measurable
functions.

\bigskip

{\bf Lemma 4.3.\/}{\it \/ Let $\Gamma$ be a Borel space with a measure and
let $\psi :\Gamma \times {\boldmat R}\arr {\boldmat R}$
      be such function that \\ (i) for every $\l\in {\boldmat R}$
the function $\psi (\g ;\l )$
is measurable on $\Gamma$; \\ (ii) $\psi (\g ;\l )$ is right-continuous and
monotonely non-increasing in the second argument for almost all $\g$. \\
For any real $\a$ put
\begin{equation} \label{15}
c_{\a}(\g)= \inf\{\l :\psi (\g ;\l )\leq \a\}.
\end{equation}
Then the function $c_{\a}(\g)$ is measurable.}

\bigskip

{\bf Proof.\/} We have to show that for any $\beta\in {\boldmat R}$
the set $V=\{ \g : c_{\a}(\g)\leq \beta\}$
must be measurable. But from the definition of $c_{\a}(\g)$ and from
({\it ii\/}) we have $V=\{ \g : \inf \{ \l :\psi (\g ;\l)\leq\a\}
\leq\beta\}=\{\g :\psi (\g ;\beta )\leq\a\}$.
By ({\it i\/}) we are done.$\quad\bullet$

\bigskip

Recall that the operator $K$ is decomposable over $\Gamma$. Let
\[
  P_{1}(\g ;\l )=\chi_{(\l  ;+\i )}(K(\g));
\]
\[
  P_{2}(\g ;\l )=\chi_{[\l ;+\i )}(K(\g))
\]
be the spectral projections of the operator $K(\g)$ corresponding
to the sets $(\l ;+\i)$ and $[\l ;+\i)$ respectively. Put
\[
P_{1}(\l)=\chi_{(\l ;+\i)}(K);\qquad P_{2}(\l)=\chi_{[\l ;+\i
)}(K)
\]
and $\phi(\g ;\l)=\bar T (P_{1}(\l))$. Then this function satisfies the
conditions of the lemma 4.3, therefore the function
\begin{equation} \label{16}
\l(\g)=\inf \{\l :\phi(\g ;\l)\leq 1\}
\end{equation}
is measurable.

\bigskip

Now we want to define two new projections in $H^{*}_{A}$:
\begin{eqnarray} \label{17}
               P_{1} & = & \int^{\oplus}_{\Gamma}\chi_{(\l(\g);+\i)}
               (K(\g))\, d\g \ =\ \int^{\oplus}_{\Gamma}P_{1}
               (\g ;\l(\g))\, d\g ; \nonumber \\ & & \\
               P_{2} & = & \int^{\oplus}_{\Gamma}\chi_{[\l(\g);+\i)}
               (K(\g))\, d\g \ =\ \int^{\oplus}_{\Gamma}P_{2}
               (\g ;\l(\g))\, d\g \nonumber
\end{eqnarray}
and we have to check correctness of this definition.

\bigskip

{\bf Lemma 4.4.\/}{\it \/ The operator-valued functions $P_{1}(\g ;\l(\g))$
and $P_{2}(\g ;\l(\g))$
are measurable.}

\bigskip

{\bf Proof.\/} It is understood that the $W^{*}$-algebra $A$ is acting on the
direct integral of Hilbert spaces $H=\int^{\oplus}_{\Gamma}H(\g)\, d\g$
with the
scalar product $(\cdot ,\cdot)$. We have to show  that the function
\begin{equation} \label{18}
\g\longmapsto \Bigl(P_{1}(\g ;\l(\g))\, \xi (\g), \xi (\g)\Bigr)
\end{equation}
is measurable for all $\xi =\int^{\oplus}_{\Gamma}\xi (\g)\, d\g \in H$.
By the theorem XIII.85 of \cite{rs} the function
\[
  \psi(\g ;\l) =\Bigl(P_{1}(\g ;\l(\g))\, \xi(\g),\xi(\g)\Bigr)
\]
satisfies the conditions of the lemma 4.3. Measurability of (\ref{18})
follows from measurability of the set $U=\{ \g :\psi(\g ;\l(\g))\leq\a\}$
for every $\a$. But from the definition of the function $c_{\a}(\g)$
(\ref{15}) one can see that
\[
U=\{\g :\l(\g)\geq c_{\a}(\g)\} =\{\g :\l(\g)-c_{\a}(\g)\geq 0\}.
\]
This set is measurable because of the measurability of function
$\l(\g)-c_{\a}(\g)$. The case of the second projection $P_{2}$ can be
handled in the same way.$\quad\bullet$

\bigskip

{\bf Corollary 4.5.\/}{\it \/ The projections $P_{1}$ and $P_{2}$ (\ref{16})
are well-defined and $\bar T (P_{1})\leq 1;\ \bar T (P_{2})\geq 1;\
 P_{1}\leq P_{2}$.}

\bigskip

These two projections define the decomposition of $H^{*}_{A}$ into
three modules:
\begin{equation} \label{19}
H^{*}_{A}=H_{-}\oplus H_{0}\oplus H_{+}
\end{equation}
where $H_{+}=\Im P_{1};\ H_{0}=\Im (P_{2}-P_{1});\ H_{-}=\Ker P_{2}$.
The operator
$K$ commutes with these projections because $K(\g)$ commutes with
the projections $P_{1}(\g ;\l(\g))$ and $P_{2}(\g ;\l(\g))$ for almost all
$\g$,
so with respect to the decomposition (\ref{19}) $K$ can be written in
the form
\[
K=\left( \begin{array}{ccc}
                      K_{+} & 0 & 0 \\
                      0 & K_{0} & 0 \\
                      0 & 0 & K_{-}
         \end{array}
\right)
\]
and $K_{0}$ for almost all $\g$ is the operator of multiplication by a
scalar $\l(\g)$, hence every submodule of $H_{0}$ is invariant for $K$.
 From the corollary 4.5 we can conclude that there exists a
projection $P$ such that $P_{1}\leq P\leq P_{2}$ and $\bar T (P)=1$. Then the
operator $K$ is diagonal also with respect  to the decomposition
$H^{*}_{A}=\Im P\oplus\Ker P$:
\[
K=\left( \begin{array}{cc}
               K_{1} & 0 \\
               0 & K'
         \end{array}
  \right) .
\]
Notice that the module $\Im P$ is isomorphic to $A$ because the
projections onto them in $H^{*}_{A}$ have the same trace $\bar T $, hence
they are equivalent \cite{t}. Let $x_{1}\in H^{*}_{A}$ be a generator of the
module $\Im P,\ \la x_{1},x_{1}\ra =1$. If it is fixed then the operator
$K_{1}:\Im P\arr \Im P$ can be viewed as the operator of
multiplication by some $\l_{1}\in A;\ K_{1}x_{1}a=x_{1}\l_{1}a$ for
$a\in A,\, x_{1}a\in\Im P$. By the theorem
2.1 the module $\Ker P$ is isomorphic to
$H^{*}_{A}$ and the operator $K'$ is obviously compact on $\Ker P$ and
the lemma 4.2 holds for it. Moreover we have the operator
inequality $K_{1}=\l_{1}\geq K'$.

    Further on by induction we can find elements $x_{i}\in H^{*}_{A}$ with
$\la x_{i},x_{j}\ra =\delta_{ij}$ and operators $\l_{i}\in A$ such that
$Kx_{i}=x_{i}\l_{i}$ and $\l_{i+1}\leq \l_{i}$.
 Denote by $N$ the $A$-module generated by these
elements $x_{i}$. Obviously $N\cong H_{A}$. Notice that the operator
$K\v_{N}$
need not to be compact. It remains to show that $N^{*}=H^{*}_{A}$.

\bigskip

{\bf Lemma 4.6.\/}{\it \/ The norm of the operators $\l_{i}$ tends to zero.}

\bigskip

{\bf Proof.\/} Since the sequence $\V\l_{i}\V$ is monotonously
non-increasing it
converges to some number $b\geq 0$. Suppose that $b\neq 0$. The operators
$\l_{i}$ as well as the other objects involved can be decomposed into
direct integrals over $\Gamma$. From construction of $\l_{i}$ we can
conclude that there exist such numbers $d_{i}(\g)$ that
\begin{equation} \label{20}
\l_{i}(\g)\geq d_{i}(\g)\geq \l_{i+1}(\g)
\end{equation}
If we decompose $x_{i}$ into a direct integral coordinatewise:
$x_{i}=\int^{\oplus}_{\Gamma}x_{i}(\g)\, d\g$
                  then for almost all $\g\quad x_{i}(\g)$ are orthonormal
in $H^{*}_{A}$ and $K(\g)x_{i}(\g)=x_{i}(\g)\l_{i}(\g)$. Define a function
$b(\g)$ as the limit of the norms $\V\l_{i}(\g)\V$ taken in $A(\g)$. We
have
\begin{equation} \label{21}
\V\l_{i}(\g)\V =\V\la K(\g)x_{i}(\g),x_{i}(\g)\ra\V\geq b(\g),
\end{equation}
where the inner product is also taken in the $A(\g)$-modules
$H^{*}_{A(\g)}$.
Let now $x$ be an element of $N$. Then it can be written in the form
\[
x(\g)=\sum_{i}x_{i}(\g)a_{i}(\g) \ \ \mbox{with some}\ \
a_{i}=\int^{\oplus}_{\Gamma}a_{i}(\g)\, d\g\in A.
\]
If $\la x,x\ra =1$ then for almost all $\g$
\begin{equation} \label{22}
\sum_{i}a^{*}_{i}(\g)a_{i}(\g)=1.
\end{equation}
{}From (\ref{20}) and (\ref{21}) we can conclude that for all $i$ the operator
inequality $\l_{i}(\g)\geq b(\g)$ holds. Therefore
\begin{eqnarray*}
\la K(\g)x(\g),x(\g)\ra & = & \sum_{i}a^{*}_{i}(\g)\l_{i}(\g)a_{i}(\g) \\
&\geq & \sum_{i}a^{*}_{i}(\g)a_{i}(\g)b(\g)\  =\  b(\g)
\end{eqnarray*}
due to (\ref{22}) and $b(\g)$ being a scalar. Further on we obtain that
\[
  \V\la Kx,x\ra\V =\ess\,\sup\V\la K(\g)x(\g),x(\g)\ra\V\geq
\ess\,\sup b(\g)=b
\]
and as by supposition $b>0$, so
\begin{equation} \label{23}
\V\la Kx,x\ra\V\geq b
\end{equation}
for any $x\in N$  with $\la x,x\ra =1$.
Now consider the projection $P_{n}:N\arr L_{n}(A)$. If the spectrum
of this operator would be separated from zero then $P_{n}$ would be
an inclusion of the module $N$ into the module $L_{n}$, but it is
impossible for finite $W^{*}$-algebras. Therefore for any $\e >0$ we can
find $x\in N$ with $\la x,x\ra =1$ such that $\V P_{n}x\V <\e$. Put
$x'=P_{n}x;\ x''=x-x'$.
 We have $\V x'\V <\e ;\ \V x''\V\leq 1$. Estimate the norm of $\la
Kx,x\ra$:
\begin{eqnarray*}
\V\la Kx,x\ra\V & \leq & \V\la Kx',x'\ra\V\  +\  2\,\V\Re\la
Kx',x''\ra\V\ +\ \V\la Kx'',x''\ra\V \\
& \leq & \V K\V\,\V x'\V^{2}\ \ +\  2\,\V K\V\,\V x'\V\,
\V x''\V\  +\ \V\la Kx'',x''\ra\V \\
& \leq & \V K\V\e^{2}\  +\ 2\,\V K\V\e\ +\ \V\la Kx'',x''\ra\V .
\end{eqnarray*}
As $x''\p L_{n}$, so due to compactness of $K$ we have $\V\la Kx'',x''\ra\V
<\e$
   for $n$ great enough. Hence $\V\la Kx,x\ra\V <\e'$ where
$\e'=\V K\V\e^{2} +2\V K\V\e +\e$.
    Choosing $\e$ small enough this estimate contradicts (\ref{23}), so
our supposition $b>0$ is wrong.$\quad\bullet$

    We have proved that the norm of the restriction of $K$ to the
orthogonal complement to $x_{1},\ldots ,x_{n}$ tends to zero. It means that
if $x\in N^{\p}$ then $\V Kx\V =0$. But $\Ker K=0$, so $N^{\p}=0$ and by the
proposition 2.5 we have $N^{\p}=H^{*}_{A}$.$\quad\bullet$

\bigskip

The ``eigenvalues'' $\l_{i}$ of $K$ are obviously not uniquely
determined and the same is true for the ``eigenvectors'' $x_{i}$. If
for example we take $x'_{i}=x_{i}u_{i}$ with unitaries $u_{i}\in A$ then the
``eigenvalues'' of $K$ will be the operators
$\l'_{i}=u^{*}_{i}\l_{i}u_{i}$. The
other reason of non-uniqueness is absence of order relation even in
commutative $W^{*}$-algebras. For example if $A=L^{\i}(X)$ and if
$\l_{i}=f(x)$;$\,\l_{j}=g(x)$ are such functions that for some
$x\quad f(x)>g(x)$ and for
some other $x$ the inverse inequality holds then the functions
$\max (f(x),g(x))$ and $\min (f(x),g(x))$ are also ``eigenvalues''.
Nevertheless the next proposition shows that putting the
``eigenvalues'' in some order provides their uniqueness.

\bigskip

{\bf Proposition 4.7.\/}{\it \/ Let $\l_{i}$ and $x_{i}$ be
as constructed in the theorem
4.1, and let $\mu_{i}$ be the ``eigenvalues'' of $K$ corresponding to
another ``basis'' $\{y_{i}\}$ of $H^{*}_{A}$. If for any unitaries
$v_{i}\in A$
and for all $i$ we have $v^{*}_{i}\mu_{i}v_{i}\geq
v^{*}_{i+1}\mu_{i+1}v_{i+1}$ then $\l_{i}$ and $\mu_{i}$
coinside up to unitary equivalence.}

\bigskip

{\bf Proof.\/} One can easily check that by supposition we have
$\inf\Sp\mu_{i}(\g)\geq\sup\Sp\mu_{i+1}(\g)$ in factor $A(\g)$
for almost all $\g\in\Gamma$. So
the projections in $H^{*}_{A}$ onto the modules generated by $y_{i}$ are
spectral projections for $K$. Denote the projection onto
$\Span_{A}(y_{1})$
by $Q$. Then obviously $P_{1}\leq Q\leq P_{2}$ where $P_{1}$,
$P_{2}$ are defined by
(\ref{17}). We can decompose $Q$ into the sum $Q=P_{1}\oplus R$ and the
projection $P$ onto $\Span_{A}(x_{1})$ into the sum $P=P_{1}\oplus S$ where
$R$ and $S$ are also projections. As $\bar T (P)=\bar T (Q)=1$,
so $R$ and $S$ are
equivalent and $\Im R\cong\Im S$. This module isomorphism commutes
with the action of $K$ because the restriction of $K$ onto these
modules is scalar and coinsides with $d_{1}(\g)$ for almost all $\g
\in \Gamma$.
So there exists a unitary $u_{1}\in A$ realizing this isomorphism
between $\Im P$ and $\Im Q$ such that $\l_{1}=u^{*}_{1}\mu_{1}u_{1}$.
Acting by induction we obtain unitary equivalence of the two sets of
``eigenvalues''.$\quad\bullet$

\bigskip

In the end of this section we must say a few words about
diagonalization theorem in the case if we drop out requests about
positiveness and absence of kernel for $K$. If $K$ is any compact
operator in $H_{A}$ or in $H^{*}_{A}$ then $H^{*}_{A}$
can be decomposed into a direct sum $H^{*}_{A}=H_{-}\oplus
\Ker K\oplus H_{+}$ so that the restriction of $K$ onto $H_{+}$
(resp. $H_{-}$ ) is positive (resp. negative). We can find sets of
``eigenvectors'' independently in $H_{+}$ and in $H_{-}$ but we need to
drop out the demand for these ``eigenvectors'' to be units, i.e. the
inner squares of such vectors are some projections but not
necessarily unities.
It is shown in \cite{fm} that any compact self-adjoint operator acting in
an autodual Hilbert module over a $W^{*}$-algebra can be diagonalized,
but its ``eigenvectors'' are not units and its ``eigenvalues'' are
not uniques up to unitary equivalence.

\section{Quadratic forms on $H^{*}_{A}$ related to
                    self-adjoint operators}

\setcounter{equation}{0}

Quadratic forms play an important role in the classical
operator theory in Hilbert spaces. If $B$ is a $C^{*}$-algebra with a
faithful finite trace $\t$ and $D$ is a self-adjoint operator acting
on a Hilbert $B$-module $M$ then a quadratic form on $M$ can be
defined as $Q(x)=\t\Bigl(\la Dx,x\ra \Bigr)$ for $x\in M$. We shall see in
this section that this $C^{*}$-module quadratic form behaves itself
like a usual one. Recall that by $B_{1}(M)$ we denote the unit ball of
$M$.

\bigskip

{\bf Proposition 5.1.\/}{\it \/ Let $D$ be a positive operator in $M$ with
$\Ker D=0$
and let the quadratic form $Q(x)$ reach its supremum on $B_{1}(M)$ at
some vector $x$. Then $\la x,x\ra$ is a projection.}

\bigskip

{\bf Proof.\/} Denote $\la x,x\ra$ by $h\in A$. By definition we have
$\V h\V\leq 1$;$\,h>0$;$\,h^{*}=h$. Suppose that the spectrum of $h$ contains
some number $c$ besides zero and unity. Define a function $f(t)$ on $[0;1]
\supset\Sp h$ by
\[
f(t)=\left\{
           \begin{array}{c}
           \frac{1}{\sqrt{\e}},\quad 0\leq t\leq\e , \\
           \frac{1}{\sqrt{t}},\quad \e\leq t\leq 1,
           \end{array}
     \right.
\]
where $0<\e <c$. Put $a=f(h)$ and $x'=xa$. Then $\la x',x'\ra
=aha=ha^{2}$.
This operator is equal to the value of the function $t\cdot f^{2}(t)$
calculated for the operator $h$. As $t\cdot f^{2}(t)\leq 1$ for $0\leq
t\leq 1$,
so $\V ha^{2}\V\leq 1$ and $x'$ lies in $B_{1}(M)$. By supposition
$\la Dx,x\ra$
is a positive operator. Denote it by $k^{2}$ with $k\geq 0$. Then
\begin{eqnarray*}
Q(x')-Q(x) & = & \t\Bigl(\la Dxa,xa\ra -\la Dx,x\ra\Bigr)\ =\
\t (ak^{2}a-k^{2}) \\
& = & \t (a^{2}k^{2}-k^{2})\ = \ \t (ka^{2}k-k^{2}) \\
& = &  \t (k(a^{2}-1)k).
\end{eqnarray*}
By definition $a^{2}-1$ is also positive and we denote it by $b^{2}$ with
$b\geq 0$. Then
\[
\t (k(a^{2}-1)k)=\t (kb^{2}k)=\t \Bigl(b\la Dx,x\ra b\Bigr)=
\t \Bigl(\la Dxb,xb\ra\Bigr)
\]
and thus
\begin{equation} \label{24}
Q(x')-Q(x)=\t \Bigl(\la Dxb,xb\ra\Bigr)
\end{equation}
But $\la xb,xb\ra =b\la x,x\ra b=bhb$ and the operator $bhb$
corresponds to the function $t(f(t)-1)$. This function differs
from zero when $t=c$, therefore the operator $bhb$ differs from
zero, and it means that $xb$ also differs from zero. Notice that
the operator $D^{1/2}$ is well-defined and $\Ker D^{1/2}=0$. Therefore
\[
\t \Bigl(\la Dxb,xb\ra \Bigr)=\t \Bigl(\la D^{1/2}xb,D^{1/2}xb\ra \Bigr)>0,
\]
hence from (\ref{24}) we obtain the inequality $Q(x')-Q(x)>0$ and it
contradicts the supposition that $Q(x)$ is the supremum of the
quadratic form $Q$ on $B_{1}(M)$. So we have proved that $\Sp h$ does
not contain any other number except zero and unity, hence $h$ is a
projection.$\quad\bullet$

\bigskip

{\bf Proposition 5.2.\/}{\it \/ Let $x\in B_{1}(M)$ be a vector at
which the quadratic
form $Q$ reaches its supremum on $B_{1}$ and let $L\subset M$ be a
submodule generated by $x$. Then $M=L\oplus L^{\p}$ and $L$ and
$L^{\p}$ are
$D$-invariant submodules, i.e. $DL\subset L$; $DL^{\p}\subset L^{\p}$.}

\bigskip

{\bf Proof.\/} By the previous proposition $\la x,x\ra$ is a projection,
hence
$L$ is a projective module and by the Dupr\'e-Fillmore theorem \cite{df}
$M=L\oplus L^{\p}$. Let $y\in L^{\p}$; $\V y\V =1$. Put
$x_{t}=x\cos t +y\sin t$. As $x_{0}=x$
is a point of maximum for $Q$, so $\frac{d}{dt}Q(x_{t})=0$ when $t=0$. It
is easy to see that
\[
\frac{d}{dt}Q(x_{t})\v_{t=0}=\t \Bigl(\la Dx,y\ra +\la Dx,y\ra^{*}\Bigr).
\]
If $\la Dx,y\ra\neq 0$ then put $a=\frac{1}{\V\la Dx,y\ra\V}
\la Dx,y\ra$. Obviously $\V a\V =1$.
Put further $z=ya^{*}$ and $\bar x_{t}=x\cos t +z\sin t$. Then
\begin{eqnarray*}
0 & = & \frac{d}{dt}Q(\bar x_{t})\v_{t=0}\ =\ \t \Bigl(\la Dx,z\ra +
\la Dx,z\ra^{*}\Bigr) \\
& = & \t \Bigl(\la Dx,y\ra a^{*} +a\la Dx,y\ra^{*}\Bigr) \ =
\ \t \Bigl(2aa^{*}\cdot\V\la Dx,y\ra\V\Bigr) \\
& = & 2\,\V\la Dx,y\ra\V
\cdot\t (aa^{*}).
\end{eqnarray*}
{}From the faithfulness of $\t$ we obtain $a=0$, hence $Dx$ is
orthogonal to any $y\in L^{\p}$, so $Dx\in L$. By self-adjointness of
$D$ we have $\la Dx,y\ra =0$ for any $y\in L^{\p}$, so
$DL^{\p}\subset L^{\p}$.$\quad\bullet$

\bigskip

{\bf Proposition 5.3.\/}{\it \/ Let $x\in B_{1}(H_{B})$ be a vector
at which the
quadratic form $Q$ reaches its supremum on $B_{1}(H_{B})$. Then
$\la x,x\ra =1$.}

\bigskip

{\bf Proof.\/} If $\la x,x\ra$ is less than unity then there exists
$y\in x^{\p}$
such that $\V y\V\leq 1$, $y\neq 0$ and $yq=y$ where $q=1-\la x,x\ra$ is a
projection. Then $\la x+y,x+y\ra =\la x,x\ra +\la y,y\ra\leq 1$, so
$x+y\in B_{1}(H_{B})$.
But $Q(x+y)=Q(x)+Q(y)$ by the previous proposition and as $y\neq 0$
and $\Ker D=0$, so $Q(y)>0$, hence $Q(x+y)>Q(x)$. This
contradiction proves the proposition.$\quad\bullet$

\bigskip

We call an operator $D$ in $M$ diagonalizable if it possesses
a ``basis'' consisting of ``eigenvectors''.

\bigskip

{\bf Proposition 5.4.\/}{\it \/ Let an operator $D$ in $M$ be positive and
diagonalizable. If for its ``eigenvalues'' one has
$\Sp\l_{i}\geq\Sp\l_{i+1}$
then the supremum of the quadratic form $Q$ on $B_{1}(M)$ is reached
at the first ``eigenvector'' $x_{1}$ and is equal to $\t(\l_{1})$.}

\bigskip

{\bf Proof.\/} Any $x\in B_{1}(M)$ can be decomposed:
$x=\sum_{i}x_{i}a_{i}$ with $a_{i}\in B$.
Then
\begin{eqnarray*}
Q(x) & = & \t \Bigl(\la Dx,x\ra\Bigr)\ \ =
\ \ \t\biggl(\sum_{i}a_{i}^{*}\la Dx_{i},
x_{i}\ra a_{i}\biggr) \\
& = & \t\biggl(\sum_{i}a^{*}_{i}\l_{i}a_{i}\biggr)\ \ \leq\ \ \t(a^{*}_{1}
\l_{1}a_{1})\ +\ \sum_{i>1}\t(a^{*}_{i}\l_{i}a_{i}).
\end{eqnarray*}
Let $\Sp\l_{1}\geq d\geq\Sp\l_{2}$. Then
\[
Q(x)=\t(a^{*}_{1}\l_{1}a_{1})+\sum_{i>1}\t(a^{*}_{i}da_{i})\leq
\t(a^{*}_{1}\l_{1}a_{1})+d\t(1-a^{*}_{1}a_{1})
\]
because the inequality $\sum_{i}a^{*}_{i}\leq 1$ follows from $\V
x\V\leq 1$.
Further on
\begin{eqnarray*}
Q(x) & \leq & \t(a^{*}_{1}\l_{1}a_{1})+d(1-\t(a^{*}_{1}a_{1}))=
\t(a^{*}_{1}\l_{1}a_{1}-a^{*}_{1}da_{1})+d \\
& = & \t(a^{*}_{1}(\l_{1}-d)a_{1})+d=\t\Bigr((\l_{1}-d)^{1/2}a_{1}a^{*}_{1}
(\l_{1}-d)^{1/2}\Bigr)+d \\
& \leq & \V
a_{1}\V^{2}\cdot\t(\l_{1}-d)+d\leq\t(\l_{1}-d)+d=\t(\l_{1}).
\end{eqnarray*}
So $\t(\l_{1})$ is the supremum of $Q(x)$ on $B_{1}(M)$ and it is
reached on $x_{1}$.$\quad\bullet$

\section{Perturbated Schr\"odinger operator with
            irrational magnetic flow as an operator acting in
                           a Hilbert module}

\setcounter{equation}{0}

In this section we consider the perturbated Schr\"odinger
operator with irrational magnetic flow
\begin{equation} \label{25}
\left(i\frac{\partial}{\partial x}+2\pi\th y\right)^{2}
-\frac{\partial^{2}}{\partial y^{2}}+W(x,y)
\end{equation}
with a double-periodic perturbation $W(x,y)=W(x+1,y)=W(x,y+1)$.
This operator has been studied in a number of papers (see \cite{ly},\cite{n}).
Applying to the operator (\ref{25}) the Fourier transform in the
variable $x\ \ (x\arr\xi)$ and the change of variables:
$t=-\frac{\xi}{2\pi}+\th y$; $s=\frac{\xi}{2\pi}$
            we obtain the operator
\begin{equation} \label{26}
     D=\Delta +W
\end{equation}
with
\begin{equation} \label{27}
  \Delta =\th^{2}\left(\left(\frac{2\pi t}{\th}\right)^{2}-
     \frac{\partial^{2}}{
  \partial t^{2}}\right)
\end{equation}
and
\[
W=\sum_{k,l}w_{kl}\,T^{k}_{t}\,T^{-k}_{s}e^{2\pi ilt/\th}e^{2\pi ils/\th}
\]
where $T_{t}$ (resp. $T_{s}$) denotes the unit translation in variable
$t$
(resp. $s$ ), $T_{t}\,\phi(t,s)=\phi(t+1,s)$, and $w_{kl}$ denote the
Fourier series coefficients of the function $W(x,y)$. We suppose
that the function $W(x,y)$ is such that $\sum_{k,l}\v w_{kl}\v
<\i$.
    Let $A_{\th}$ be the $C^{*}$-algebra generated by two non-commuting
unitaries $U$ and $V$ such that $UV=e^{2\pi i\th}VU$  \cite{br},\cite{co}
and let
$A^{\i}_{\th}\subset A_{\th}$ be its ``infinitely smooth''
subalgebra of elements of
the form $\sum_{k,l}a_{kl}U^{k}V^{l}$ where coefficients $a_{kl}$ are of
rapid decay. The Schwartz space $S({\boldmat R})$ of functions of rapid decay
on ${\boldmat R}$ can be made \cite{co} a projective right $A^{\i}
_{\th}$-module with one
generator. We denote this module by $M^{\i}$. The action of
$A^{\i}_{\th}$ on $M^{\i}$
is given by formulas
\[
(\phi U)(t)=\phi(t+\th);\quad (\phi V)(t)=e^{2\pi it}\phi(t)
\]
for $\phi(t)\in M^{\i}$. The module $M^{\i}$ is generated by a projection
$p\in A^{\i}_{\th}$; $M^{\i}\cong pA^{\i}_{\th}$ with $\t(p)=\th$
and as $M^{\i}\subset A_{\th}$ so $M^{\i}$
inherits the norm from $A_{\th}$. Its closure $M=M^{\i}\otimes_{
A^{\i}_{\th}}A_{\th}$ in
this norm is a Hilbert $A_{\th}$-module. Notice that there exists in
$S({\boldmat R})\subset L^{2}({\boldmat R})$ the orthonormal basis
$\{\phi_{i}(t)\}$ consisting of the
eigenfunctions of the operator $\Delta$ (\ref{27}), and the functions from
$S({\boldmat R}^{2})=M^{\i}\mathbin{\hat{\otimes}} M^{\i}$ can be represented
as series $\sum_{i}\phi_{i}(t)m_{i}(s)$
with $m_{i}(s)\in M^{\i}$. Define the $A_{\th}$-valued inner product on
$S({\boldmat R}^{2})$ by formula
\[
\Bigl <\sum_{i}\phi_{i}(t)m_{i}(s),\sum_{j}\phi_{j}(t)n_{j}(s)
\Bigr > =\sum_{i}\la m_{i}(s),n_{i}(s)\ra
\]
where $n_{j}(s)\in M^{\i}$. By $S({\boldmat R};M)$ (resp.
$L^{2}({\boldmat R};M)$ ) we denote
the Schwartz space of functions (resp. the space of
square-integrable functions) with the values in the Banach space $M$.
The inclusion
\[
S({\boldmat R}\times{\boldmat R})\hookrightarrow
S({\boldmat R};M)\hookrightarrow L^{2}({\boldmat R};M)\cong N
\]
allows us to consider $S({\boldmat R}^{2})$ as a dense subspace in the Hilbert
module
\[
N=\{(m_{i}):\sum_{i}\la m_{i},m_{i}\ra\ \mbox{converges in}
\ A_{\th}\}
\]
(this module is
often denoted by $l_{2}(M)$ ). One can see that the module $M$ is
full, i.e. $\la M,M\ra =A_{\th}$ because the $C^{*}$-algebra $A_{\th}$ is
simple
and $\la M,M\ra$ must be its ideal. By the results of \cite{df} one has
$N\cong H_{A_{\th}}$.

\bigskip

{\bf Theorem 6.1.\/}{\it \/ The operator $D$ (\ref{26}) is a self-adjoint
unbounded operator in $N$ with a dense domain.}

\bigskip

{\bf Proof \/} consists of the five following steps.

\medskip   \noindent
1. Let $N_{1}\subset N$ be a subspace of sequences $(m_{i})$ such that the
series $\sum_{i}i^{2}\la m_{i},m_{i}\ra$ converges in norm to an element
of $A_{\th}$.
If $\xi\in N_{1}$; $\xi =\sum_{i}\phi_{i}(t)m_{i}$ then
$\Delta\xi =\sum_{i}(2i-1)\th\phi_{i}(t)m_{i}$
and the series $\sum_{i}(2i-1)^{2}\th^{2}\la m_{i},m_{i}\ra$ converges in
$A_{\th}$,
therefore $\Delta$ is an unbounded operator in $N$ with the dense domain
$N_{1}$.

\medskip \noindent
2. Here we show that the action of the operator
$C_{kl}=T^{-k}_{s}e^{2\pi ils/\th}$
can be prolonged from $M^{\i}$ to $M$. Since this operator commutes
with the action of the algebra $A^{\i}_{\th}$ on the
$A^{\i}_{\th}$-module $M^{\i}$ we
have $C_{kl}\in \End_{A^{\i}_{\th}}M^{\i}$. The image of the generator
$p$ of $M^{\i}$
can be written in the form $C_{kl}(p)=pa_{kl}\in M^{\i}$  for some
$a_{kl}\in A^{\i}_{\th}$
and we have
\[
C_{kl}(p)=C_{kl}(p^{2})=C_{kl}(p)p=pa_{kl}p.
\]
Obviously the map $m\longmapsto pa_{kl}pm=C_{kl}(m)$ can be continuously
prolonged from $M^{\i}$ to $M$. Besides that since $C_{kl}$ is a
unitary operator, we have $\V C_{kl}\V =1$ and $\V a_{kl}\V =1$.

\medskip \noindent
3. Consider now the operator $B_{kl}=T^{k}_{t}e^{2\pi ilt/\th}
\cdot C_{kl}$. It
is obviously continuous in $S({\boldmat R}^{2})$. Let $\a_{ij}$ be matrix
coefficients of decomposition with respect to the basis $\{\phi_{j}\}$ for
the operator $T^{k}_{t}e^{2\pi ilt/\th}$:
\[
T^{k}_{t}e^{2\pi ilt/\th}\phi_{i}(t)=\sum_{j}\a_{ij}\phi_{j}(t).
\]
As this operator is unitary, so
$\sum_{j}\bar\a_{ij}\a_{nj}=\delta_{in}$. Let
$\xi =\sum_{i}\phi_{i}(t)m_{i}\in N$.
Then $B_{kl}(\xi)=\sum_{i,j}\a_{ij}\phi_{j}(t)C_{kl}(m_{i})$.
Estimate its norm:
\begin{eqnarray*}
\la B_{kl}(\xi),B_{kl}(\xi)\ra & = & \sum_{i,j}\Bigl <\sum_{i}\a_{ij}
C_{kl}(m_{i}),\sum_{n}\a_{nj}C_{kl}(m_{n})\Bigr > \\
& = & \sum_{i,n,j}\bar\a_{ij}\a_{nj}\la C_{kl}(m_{i}),C_{kl}(m_{n})\ra
\\
& = & \sum_{i,n}\biggl(\sum_{j}\bar\a_{ij}\a_{nj}\biggr)
\la C_{kl}(m_{i}),C_{kl}(m_{n})\ra \\
& = & \sum_{i,n}\delta_{in}\la C_{kl}(m_{i}),C_{kl}(m_{n})\ra \\
& = & \sum_{i}\la C_{kl}(m_{i}),C_{kl}(m_{i})\ra \\
& = & \sum_{i}(a_{kl}m_{i})^{*}a_{kl}m_{i} = \sum_{i}m^{*}_{i}
a_{kl}^{*}a_{kl}m_{i} \\
& \leq & \V a_{kl}\V^{2}\sum_{i}m_{i}^{*}m_{i} =\sum_{i}
m_{i}^{*}m_{i} = \la \xi ,\xi\ra .
\end{eqnarray*}
Hence $\V B_{kl}\V\leq 1$ and it is a continuous operator in $N$.

\medskip  \noindent
4. We have
\[
\V W\V\leq\sum_{k,l}\V w_{kl}B_{kl}\V\leq\sum_{k,l}\v w_{kl}\v\cdot
\V B_{kl}\V\leq\sum_{k,l}\v w_{kl}\v .
\]
By our supposition the last sum is finite, hence $W$ is continuous in $N$.

\medskip \noindent
5. It remains to show that $D$ commutes with the action of the
$C^{*}$-algebra $A_{\th}$ on $N$. It is obvious for the operators
$\Delta$ and $B_{kl}$.
 As the series $W=\sum_{k,l}w_{kl}B_{kl}$ converges, so $W$ also
commutes with the action of $A_{\th}$. $D$ is self-adjoint if the
function $W(x,y)$ is real-valued. $\quad\bullet$

\bigskip

Let now $A$ be a type $\mbox{II}_{1}$ factor containing $A_{\th}$ as
a weakly
dense subalgebra (cf. \cite{br}). This inclusion induces the inclusion of
$H_{A_{\th}}$ into $H_{A}$ and operators acting in $H_{A_{\th}}$ can be
prolonged to operators acting in $H_{A}$. Notice that if $\V W\V <c$
then the operator $D+c$ is invertible and its inverse $(\Delta
+W+c)^{-1}=(1+\Delta^{-1}(W+c))^{-1}\Delta^{-1}$
is compact because the operator $\Delta^{-1}$ is
compact. So by the theorem 4.1 it is diagonalizable in $H^{*}_{A}$, hence
the same is true for the operator $D$. Slightly changing the proof
of that theorem (namely taking $\th$ instead of $1$ in (\ref{16})) we can
obtain the set of ``eigenvectors'' $\{x_{i}\}$ for $D$ with
$\la x_{i},x_{i}\ra =p$.
In that case the corresponding ``eigenvalues'' $\l_{i}$ can be viewed as
elements from $\End^{*}_{A}({\cal N})$ where ${\cal N}=pA
=N\otimes_{A_{\th}}A$.

\bigskip

{\bf Problem 6.2.\/} Can the ``eigenvalues'' $\l_{i}$ be taken from the lesser
algebra $\End^{*}_{A_{\th}}(M)$
instead of $\End^{*}_{A}({\cal N})$\,?
Do these ``eigenvalues'' possess properties
resembling analyticity as they do in the commutative case when $\th$
is integer \cite{n},\cite{rs}\,?

\bigskip

If $\V W\V <\th$ then the spectrum of $D$ lies in
$\cup_{i}(2\th (i-1);2\th i)$,
therefore the spectral projections $P_{i}=P_{(2\th (i-1);\,2\th i)}
(D)$ lie in $\End^{*}_{A_{\th}}(N)$, hence the ``eigenvalues'' $\l_{i}$
of $D$ lie in $\End^{*}_{A_{\th}}(M)$.
 It was shown in \cite{m1} by the methods of perturbation
theory that if the norm of $W$ is small enough then the images of
$P_{i}$
contain ``eigenvectors'' which form a basis of $N$, hence the
operator $D$ is diagonalizable inside the module $N$.

\bigskip

{\bf Added in proof:} \, The results of \cite{st} allow us to give positive
answer to the first question of the problem 6.2.

%\section*{}
\bigskip\bigskip
\bigskip\bigskip
V.M.Manuilov \\*
Moscow State Public University \\*
Russia - 129278 Moscow \\*
P.Korchagina str.\ 22 \\*
E-mail:manuilov@math.math.msu.su

\end{document}